\documentclass[a4paper,11pt]{amsart}
\usepackage{amssymb,amscd}
\usepackage[cp1251]{inputenc}
\usepackage[T2A]{fontenc}
\theoremstyle{plain}

\newtheorem{theo}{Theorem}[section]
\newtheorem{prop}[theo]{Proposition}

\theoremstyle{definition}

\newtheorem{exam}[theo]{Example}
\newtheorem{prob}[theo]{Problem}

\theoremstyle{remark}

\numberwithin{equation}{section}

\newcommand{\field}[1]{\mathbb{#1}}
\def\C{\field{C}}
\newcommand{\R}{\field{R}}

\newcommand{\Z}{\field{Z}}

\def\D{\mathbb D}
\def\G{G}
\def\K{K}
\def\Sph{\mathbb S}
\def\T{T}

\DeclareMathOperator{\Tor}{Tor}

\def\le{\leqslant}
\def\ge{\geqslant}

\def\dbs{/\!\!/}
\def\cone{\mathop{\mathrm{cone}}}
\def\kn{\mathop{\mbox{\it KN\/}}\nolimits}

\begin{document}

\title{Topology of Kempf--Ness sets for algebraic torus actions}

\author{Taras Panov}\thanks{The work was supported by the Russian
Foundation for Basic Research (grants no.~08-01-00541,
08-01-91855-KO), the State Programme for the Support of Leading
Scientific Schools (grant 1824.2008.1), and P.~Deligne's 2004
Balzan prize in mathematics.}
\address{Department of Geometry and Topology, Faculty of Mathematics and Mechanics, Moscow State University,
Leninskiye Gory, Moscow 119992, Russia\newline \emph{and}
\newline Institute for Theoretical and Experimental Physics, Moscow
117259, Russia}
\email{tpanov@mech.math.msu.su}


\begin{abstract}
In the theory of algebraic group actions on affine varieties, the
concept of a Kempf--Ness set is used to replace the categorical
quotient by the quotient with respect to a maximal compact
subgroup. By making use of the recent achievements of ``toric
topology'' we show that an appropriate notion of a Kempf--Ness set
exists for a class of algebraic torus actions on quasiaffine
varieties (coordinate subspace arrangement complements) arising in
the Batyrev--Cox ``geometric invariant theory'' approach to toric
varieties. We proceed by studying the cohomology of these
``toric'' Kempf--Ness sets. In the case of projective non-singular
toric varieties the Kempf--Ness sets can be described as complete
intersections of real quadrics in a complex space.
\end{abstract}

\maketitle

\section{Introduction}
The concept of a Kempf--Ness set plays an important role in the
geometric invariant theory, as explained, for example,
in~\cite[\S6.12]{po-vi89} or~\cite{schw89}. Given an affine
variety $S$ over $\C$ with an action of a reductive group $\G$,
one can find a compact subset $\kn\subset S$ such that the
categorical quotient $S\dbs\G$ is homeomorphic to the quotient
$\kn\!\!/\K$ of $\kn$ by a maximal compact subgroup $\K\subset\G$.
Another important property of the Kempf--Ness set $\kn$ is that it
is a $\K$-equivariant deformation retract of $S$.

Our aim here is to extend the notion of a Kempf--Ness set to a
class of algebraic torus actions on complex quasiaffine varieties
(coordinate subspace arrangement complements) arising in the
theory of toric varieties. Although our Kempf--Ness sets cannot be
defined exactly in the same way as in the affine case, they
possess the two above characteristic properties. In the case of a
projective toric variety, our Kempf--Ness set can be identified
with the level surface for the moment map corresponding to a
compact torus action on the complex space~\cite[\S4]{cox97}. The
toric Kempf--Ness sets also constitute a particular subclass of
\emph{moment-angle complexes}~\cite{bu-pa02}, which opens new
links between toric topology and geometric invariant theory.

In Section~\ref{knsav} we review the notion of Kempf--Ness sets
for reductive groups acting on affine varieties. In
Section~\ref{altoa} we outline the ``geometric invariant theory''
approach to toric varieties as quotients of algebraic torus
actions on coordinate subspace arrangement complements, and
introduce a toric Kempf--Ness set using our construction of
moment-angle complexes. In Section~\ref{norfa} we restrict our
attention to torus actions arising from normal fans of convex
polytopes. In this case the corresponding Kempf--Ness set admits a
transparent geometric interpretation as a complete intersection of
real quadratic hypersurfaces. The quotient toric variety is
projective, and the Kempf--Ness set is the level surface for an
appropriate moment map, thereby extending the analogy with the
affine case even further in Section~\ref{prvar}. In the last
Section~\ref{cohom} we give a description of the cohomology ring
of the Kempf--Ness set. As it is clear from an example provided,
our Kempf--Ness sets may be quite complicated topologically; many
interesting phenomena occur even for the torus actions
corresponding to simple 3-dimensional fans.

I am grateful to Mikiya Masuda for his great hospitality during my
stay in Osaka and the insight gained from numerous discussions on
toric manifolds. My thanks also go to Gerald Schwarz for
introducing me to the notion of a Kempf--Ness set and several
other important concepts of Geometric Invariant Theory. I
gratefully acknowledge many important remarks and comments
suggested by Ivan Arzhantsev, who read the first version of the
manuscript. Finally, I thank the referee whose valuable comments
and corrections improved the text significantly. Part of this work
was presented at the conference on Contemporary Geometry and
Related Topics in Belgrade (June, 2005), and I take this
opportunity to thank the organisers for their hospitality.

\section{Kempf--Ness sets for affine varieties}\label{knsav}
We start by briefly reviewing quotients and Kempf--Ness sets of
reductive group actions on affine varieties. The details can be
found in~\cite[\S6.12]{po-vi89} and~\cite{schw89}.

Let $\G$ be a \emph{reductive} algebraic group acting on an affine
variety~$X$. As $G$ is non-compact, taking the standard (or
\emph{geometric}) quotient with the quotient topology may result
in a badly behaving space (e.g. it may fail to be Hausdorff). An
alternative notion of a \emph{categorical quotient} remedies this
difficulty and ensures that the result always lies within the
category of algebraic varieties.

Let $\C[X]$ be the algebra of regular functions on~$X$, so that
$X=\mathop{\mathrm{Spec}}\C[X]$. Denote $X\dbs\G$ the complex
affine variety corresponding to the subalgebra $\C[X]^\G$ of
$\G$-invariant polynomial functions on $X$, and let $\varrho\colon
X\to X\dbs G$ be the morphism dual to the inclusion $\C[X]^\G\to
\C[X]$. Then $\varrho$ is surjective and establishes a bijection
between closed $\G$-orbits of $X$ and points of $X\dbs \G$.
Moreover, $\varrho$ is universal in the class of morphisms from
$X$ constant on $\G$-orbits in the category of algebraic varieties
(which explains the term ``categorical quotient''). The
categorical quotient coincides with the geometric one if and only
if all $G$-orbits are closed.

\begin{exam}
Consider the standard $\C^*$-action on $\C$ (here $\C^*$ is the
multiplicative group of complex numbers). The categorical quotient
$\C\dbs\C^*$ is a point, while $\C/\C^*$ is a non-Hausdorff
two-point space.
\end{exam}

Let $\rho\colon \G\to\mathop{\mathrm{GL}}(W)$ be a representation
of $\G$, let $\K$ be a maximal compact subgroup of $\G$, and let
$\langle\:,\,\rangle$ be a $\K$-invariant hermitian form on $W$
with associated norm $\|\:\|$. Given $v\in W$, consider the
function $F_v\colon \G\to\R$ sending $g$ to $\frac12\|gv\|^2$. It
has a critical point if and only if the orbit $\G v$ is closed,
and all critical points of $F_v$ are
minima~\cite[Th.~6.18]{po-vi89}. Define the subset $\kn\subset W$
by one of the following equivalent conditions:
\begin{align}\label{eqcns}\notag
  \kn&=\{v\in W\colon (dF_v)_e=0\}\qquad\text{($e\in\G$ is the unit)}\\\notag
   &=\{v\in W\colon T_v \G v\perp v\}\\\notag
   &=\{v\in W\colon\langle\gamma v,v\rangle=0\text{ for all }\gamma\in\mathfrak g\}\\
   &=\{v\in W\colon\langle\kappa v,v\rangle=0\text{ for all }\kappa\in\mathfrak k\},
\end{align}
where $\mathfrak g$ (resp. $\mathfrak k$) is the Lie algebra of
$\G$ (resp. $\K$) and we consider $\mathfrak k\subset\mathfrak
g\subset\mathop{\mathrm{End}}(W)$. Therefore, any point $v\in \kn$
is a closest point to the origin in its orbit $\G v$. Then $\kn$
is called the \emph{Kempf--Ness set} of~$W$.

We may assume that the affine $\G$-variety $X$ is equivariantly
embedded as a closed subvariety in a $\G$-representation space
$W$. Then the \emph{Kempf--Ness set} $\kn_X$ of $X$ is defined as
$\kn\cap\:X$.

The importance of Kempf--Ness sets for the study of orbit
quotients is due to the following result, whose proof can be found
in \cite[(4.7), (5.1)]{schw89}.

\begin{theo}\label{knspr}
{\rm(a)} The composition $\kn_X\hookrightarrow X\to X\dbs \G$ is
proper and induces a homeomorphism $\kn_X/\K\to X\dbs \G$.

{\rm(b)} There is a $\K$-equivariant deformation retraction of $X$
to $\kn_X$.
\end{theo}

\section{Algebraic torus actions}\label{altoa}
Let $N\cong\Z^n$ be an integral lattice of rank $n$, and
$N_\R=N\otimes_\Z\R$ the ambient real vector space. A convex
subset $\sigma\in N_\R$ is called a \emph{cone} if there exist
vectors $a_1,\ldots,a_k\in N$ such that
\[
  \sigma=\{\mu_1a_1+\ldots+\mu_ka_k\colon\mu_i\in\R,\;\mu_i\ge0\}.
\]
If the set $\{a_1,\ldots,a_k\}$ is minimal, then it is called the
\emph{generator set} of~$\sigma$. A cone is called \emph{strongly
convex} if it contains no line; all the cones below are assumed to
be strongly convex. A cone $\sigma$ is called \emph{regular}
(resp. \emph{simplicial}) if $a_1,\ldots,a_k$ can be chosen to
form a subset of a $\Z$-basis of $N$ (resp. an $\R$-basis of
$N_\R$). A \emph{face} of a cone $\sigma$ is the intersection
$\sigma\cap H$ with a hyperplane $H$ for which the whole $\sigma$
is contained in one of the two closed halfspaces determined by
$H$; a face of a cone is again a cone. Every generator of $\sigma$
spans a one-dimensional face, and every face of $\sigma$ is
spanned by a subset of the generator set.

A finite collection $\Sigma=\{\sigma_1,\ldots,\sigma_s\}$ of cones
in $N_\R$ is called a \emph{fan} if a face of every cone in
$\Sigma$ belongs to $\Sigma$ and the intersection of any two cones
in $\Sigma$ is a face of each. A fan $\Sigma$ is called
\emph{regular} (resp. \emph{simplicial}) if every cone in $\Sigma$
is regular (resp. simplicial). A fan
$\Sigma=\{\sigma_1,\ldots,\sigma_s\}$ is called \emph{complete} if
$N_\R=\sigma_1\cup\ldots\cup\sigma_s$.

Let $\C^*=\C\setminus\{0\}$ be the multiplicative group of complex
numbers, and $\Sph^1$ be the subgroup of complex numbers of
absolute value one. The \emph{algebraic torus}
$\T_\C=N\otimes_\Z\C^*\cong (\C^*)^n$ is a commutative complex
algebraic group with a maximal compact subgroup
$\T=N\otimes_\Z\Sph^1\cong (\Sph^1)^n$, the (compact)
\emph{torus}. A \emph{toric variety} is a normal algebraic variety
$X$ containing the algebraic torus $\T_\C$ as a Zariski open
subset in such a way that the natural action of $\T_\C$ on itself
extends to an action on~$X$.

There is a classical construction (cf.~\cite{dani78}) establishing
a one-to-one correspondence between fans in $N_\R$ and complex
$n$-dimensional toric varieties. Regular fans correspond to
non-singular varieties, while complete fans give rise to compact
ones. Below we review another construction of toric varieties as
certain algebraic quotients; it is due to several authors
(cf.~\cite{baty93}, \cite{cox95}).

In the rest of this section we assume that the one-dimensional
cones of $\Sigma$ span $N_\R$ as a vector space (this holds, e.g.,
if $\Sigma$ is a complete fan). Assume that $\Sigma$ has $m$
one-dimensional cones. We order them arbitrarily and consider the
map $\Z^m\to N$ sending the $i$th generator of $\Z^m$ to the
integer primitive vector $a_i$ generating the $i$th
one-dimensional cone. The corresponding map of the algebraic tori
fits into an exact sequence
\begin{equation}\label{ggrou}
  1\longrightarrow \G\longrightarrow (\C^*)^m\longrightarrow \T_\C
  \longrightarrow1,
\end{equation}
where $\G$ is isomorphic to a product of $(\C^*)^{m-n}$ and a
finite group. If $\Sigma$ is a regular fan and has at least one
$n$-dimensional cone, then $\G\cong(\C^*)^{m-n}$. We also have an
exact sequence of the corresponding maximal compact subgroups:
\begin{equation}\label{kgrou}
  1\longrightarrow \K\longrightarrow \mathbb T^m\longrightarrow \T
  \longrightarrow1
\end{equation}
(here and below we denote $\mathbb T^m=(\Sph^1)^m$).

We say that a subset $\{i_1,\ldots,i_k\}\subset
[m]=\{1,\ldots,m\}$ is a \emph{g-subset} if
$\{a_{i_1},\ldots,a_{i_k}\}$ is a subset of the generator set of a
cone in $\Sigma$. The collection of $g$-subsets is closed with
respect to the inclusion, and therefore forms an (abstract)
simplicial complex on the set $[m]$, which we denote $\mathcal
K_\Sigma$. Note that if $\Sigma$ is a complete simplicial fan,
then $\mathcal K_\Sigma$ is a triangulation of an
$(n-1)$-dimensional sphere. Given a cone $\sigma\in\Sigma$, we
denote by $g(\sigma)\subseteq[m]$ the set of its generators. Now
set
\[
  A(\Sigma)=\bigcup_{\{i_1,\ldots,i_k\}\text{ is not a $g$-subset}}
  \{z\in\C^m\colon z_{i_1}=\ldots=z_{i_k}=0\}
\]\\[-10pt]
and\\[-20pt]
\[
  U(\Sigma)=\C^m\setminus A(\Sigma).
\]
Both sets depend only on the combinatorial structure of the
simplicial complex $\mathcal K_\Sigma$; the set $U(\Sigma)$
coincides with the \emph{coordinate subspace arrangement
complement} $U(\mathcal K_\Sigma)$ considered
in~\cite[\S8.2]{bu-pa02}.

The set $A(\Sigma)$ is an affine variety, while its complement
$U(\Sigma)$ admits a simple affine cover, as described in the
following statement.

\begin{prop}\label{usigm}
Given a cone $\sigma\in\Sigma$, set $z^{\hat\sigma}=\prod_{j\notin
g(\sigma)}z_j$ and define
\[
  V(\Sigma)=\{z\in\C^m\colon z^{\hat\sigma}=0\text{ for all }\;\sigma\in\Sigma\}
\]\\[-15pt]
and\\[-15pt]
\[
  U(\sigma)=\{z\in\C^m\colon z_j\ne0\;\text{ if }\;j\notin g(\sigma)\}.
\]
Then $A(\Sigma)=V(\Sigma)$ and
\[
  U(\Sigma)=\C^m\setminus
  V(\Sigma)=\bigcup_{\sigma\in\Sigma}U(\sigma).
\]
\end{prop}
\begin{proof}We have
\[
  \C^m\setminus V(\Sigma)=\bigcup_{\sigma\in\Sigma}
  \{z\in\C^m\colon z^{\hat\sigma}\ne0\}=\bigcup_{\sigma\in\Sigma}U(\sigma).
\]
On the other hand, given a point $z\in\C^m$, denote by
$\omega(z)\subseteq[m]$ the set of its zero coordinates. Then
$z\in\C^m\setminus A(\Sigma)$ if and only if $\omega(z)$ is a
$g$-subset. This is equivalent to saying that $z\in U(\sigma)$ for
some $\sigma\in\Sigma$. Therefore, $\C^m\setminus
A(\Sigma)=\cup_{\sigma\in\Sigma}U(\sigma)$, thus proving the
statement.
\end{proof}

The complement $U(\Sigma)$ is invariant with respect to the
$(\C^*)^m$-action, and it is easy to see that the subgroup $\G$
from~\eqref{ggrou} acts on $U(\Sigma)$ with finite isotropy
subgroups if $\Sigma$ is simplicial (or even freely if $\Sigma$ is
a regular fan). The corresponding quotient is identified with the
toric variety $X_\Sigma$ determined by $\Sigma$. The more precise
statement is as follows.

\begin{theo}[{cf.~\cite[Th.~2.1]{cox95}}]\label{coxth}
Assume that the one-dimensional cones of $\Sigma$ span $N_\R$ as a
vector space.

{\rm(a)} The toric variety $X_\Sigma$ is naturally isomorphic to
the categorical quotient of $U(\Sigma)$ by $\G$.

{\rm(b)} $X_\Sigma$ is the geometric quotient of $U(\Sigma)$ by
$\G$ if and only if $\Sigma$ is simplicial.
\end{theo}

Therefore, if $\Sigma$ is a simplicial (in particular, regular)
fan satisfying the assumption of Theorem~\ref{coxth}, then all the
orbits of the $\G$-action on $U(\Sigma)$ are closed and the
categorical quotient $U(\Sigma)\dbs \G$ can be identified with
$U(\Sigma)/\G$. However, the analysis of the previous section does
not apply here, as $U(\Sigma)$ is \emph{not} an affine variety in
$\C^m$ (it is only quasiaffine in general). For example, if
$\Sigma$ is a complete fan, then the $\G$-action on the whole
$\C^m$ has only one closed orbit, the origin, and the quotient
$\C^m\dbs \G$ consists of a single point. In the rest of the paper
we show that an appropriate notion of the Kempf--Ness set exists
for this class of torus actions, and study some of its most
important topological properties.

Consider the unit polydisc
\[
  (\D^2)^m=\{z\in\C^m\colon|z_j|\le1\text{ for all }j\}.
\]
Given $\sigma\in\Sigma$, define
\[
  \mathcal Z(\sigma)=\{z\in(\D^2)^m\colon |z_j|=1\text{ if }j\notin
  g(\sigma)\},
\]\\[-15pt]
and\\[-20pt]
\[
  \mathcal Z(\Sigma)=\bigcup_{\sigma\in\Sigma}\mathcal Z(\sigma).
\]
The subset $\mathcal Z(\Sigma)\subseteq(\D^2)^m$ is invariant with
respect to the $\mathbb T^m$-action. (We have $\mathcal
Z(\Sigma)=\mathcal Z_{\mathcal K_\Sigma}$, where $\mathcal
Z_{\mathcal K}$ is the \emph{moment-angle complex} associated with
a simplicial complex $\mathcal K$ in~\cite[\S6.2]{bu-pa02}.) Note
that $\mathcal Z(\sigma)\subset U(\sigma)$, and therefore,
$\mathcal Z(\Sigma)\subset U(\Sigma)$ by Proposition~\ref{usigm}.

\begin{prop}
Assume that $\Sigma$ is a complete simplicial fan. Then $\mathcal
Z(\Sigma)$ is a compact $(m+n)$-manifold with a $\mathbb
T^m$-action.
\end{prop}
\begin{proof}
As $\mathcal K_\Sigma$ is a triangulation of an
$(n-1)$-dimensional sphere, the result follows
from~\cite[Lemma~6.13]{bu-pa02} (or~\cite[Lemma~3.3]{pano05}).
\end{proof}

\begin{theo}\label{zspro}
Assume that $\Sigma$ is a simplicial fan.

{\rm(a)} If $\Sigma$ is complete, then the composition $\mathcal
Z(\Sigma)\hookrightarrow U(\Sigma)\to U(\Sigma)/\G$ induces a
homeomorphism $\mathcal Z(\Sigma)/\K\to U(\Sigma)/\G$.

{\rm(b)} There is a $\mathbb T^m$-equivariant deformation
retraction of $U(\Sigma)$ to $\mathcal Z(\Sigma)$.
\end{theo}
\begin{proof}
Denote by $\cone \mathcal K_\Sigma'$ the cone over the barycentric
subdivision of $\mathcal K_\Sigma$ and by $C(\Sigma)$ the
topological space $|\cone \mathcal K_\Sigma'|$ with the dual
\emph{face decomposition}, see~\cite[\S3.1]{pano05} for details.
(If $\Sigma$ is a complete fan, then $\mathcal K_\Sigma$ is a
sphere triangulation, $C(\Sigma)$ can be identified with the unit
ball in $N_\R$, and the face decomposition of its boundary is
Poincar\'e dual to $\mathcal K_\Sigma$.) The space $C(\Sigma)$ has
a face $C(\sigma)$ of dimension $n-g(\sigma)$ for each cone
$\sigma\in\Sigma$. Set
\[
  T(\sigma)=\{(t_1,\ldots,t_m)\in \mathbb T^m\colon t_j=1\text{ if }j\notin
  g(\sigma)\}.
\]
This is a $g(\sigma)$-dimensional coordinate subgroup in $\mathbb
T^m$. As detailed in~\cite{da-ja91} and~\cite[\S3.1]{pano05}, the
set $\mathcal Z(\Sigma)$ can be described as the identification
space
\[
  \mathcal Z(\Sigma)=\bigl(\mathbb T^m\times C(\Sigma)\bigr)/{\sim},
\]
where $(t,x)\in \mathbb T^m\times C(\Sigma)$ is identified with
$(s,x)\in \mathbb T^m\times C(\Sigma)$ if $x\in C(\sigma)$ and
$t^{-1}s\in T(\sigma)$ for some $\sigma\in\Sigma$. The map of tori
$\mathbb T^m\to\T$ with kernel $\K$ induces a map of the
identification spaces
\[
  (\mathbb T^m\times C(\Sigma))/{\sim}\to(\T\times C(\Sigma))/{\sim}.
\]
Now, according to~\cite{da-ja91}, if $\Sigma$ is a complete
simplicial fan, then the latter identification space is
homeomorphic to the toric variety $X_\Sigma=U(\Sigma)/\G$. This
proves~(a). (Note that if $\Sigma$ is a regular fan, then $\K\cong
\mathbb T^{m-n}$ and the projection $\mathcal Z_\Sigma\to
X_\Sigma$ is a principal $\K$-bundle.)

Statement (b) is proved in~\cite[Th.~8.9]{bu-pa02}.
\end{proof}

By comparing this result with Theorem~\ref{knspr}, we see that
$\mathcal Z(\Sigma)$ has the same properties with respect to the
$\G$-action on $U(\Sigma)$ as the set $\kn_S$ with respect to a
reductive group action on an affine variety $S$. We therefore
refer to $\mathcal Z(\Sigma)$ as the \emph{Kempf--Ness set} of
$U(\Sigma)$.

\begin{exam}
Let $n=2$ and $e_1,e_2$ be a basis in $N_\R$.

1. Consider a complete fan $\Sigma$ having the following three
2-dimensional cones: the first is spanned by $e_1$ and $e_2$, the
second spanned by $e_2$ and $-e_1-e_2$, and the third spanned by
$-e_1-e_2$ and $e_1$. The simplicial complex $\mathcal K_\Sigma$
is a complete graph on 3 vertices (or the boundary of a triangle).
We have
\[
  U(\Sigma)=\C^3\setminus\{z\colon
  z_1=z_2=z_3=0\}=\C^3\setminus\{0\}
\]\\[-20pt]
and\\[-10pt]
\[
  \mathcal Z(\Sigma)=(\D^2\times\D^2\times\mathbb S^1)\cup
  (\D^2\times \mathbb S^1\times\D^2)\cup(\mathbb S^1\times\D^2\times\D^2)=
  \partial((\D^2)^3)\cong \mathbb S^5.
\]
The subgroup $\G$ from exact sequence~\eqref{ggrou} is the
diagonal 1-dimensional subtorus in $(\C^*)^3$, and $\K$ is the
diagonal subcircle in $\mathbb T^3$. Therefore, we have
$X_\Sigma=U(\Sigma)/\G=\mathcal Z(\Sigma)/\K=\C P^2$, the complex
projective 2-plane.

2. Now consider the fan $\Sigma$ consisting of three 1-dimensional
cones generated by vectors $e_1$, $e_2$ and $-e_1-e_2$. This fan
is not complete, but its 1-dimensional cones span $N_\R$ as a
vector space. So Theorem~\ref{coxth} applies, but
Theorem~\ref{zspro}~(a) does not. The simplicial complex $\mathcal
K_\Sigma$ consists of 3 disjoint points. The space $U(\Sigma)$ is
the complement to 3 coordinate lines in $\C^3$:
\[
  U(\Sigma)=\C^3\setminus\{z\colon
  z_1=z_2=0,\quad z_1=z_3=0,\quad z_2=z_3=0\},
\]\\[-20pt]
and\\[-15pt]
\[
  \mathcal Z(\Sigma)=(\D^2\times\Sph^1\times\Sph^1)\cup
  (\Sph^1\times\D^2\times\Sph^1)\cup(\Sph^1\times\Sph^1\times\D^2).
\]
Both spaces are homotopy equivalent to
$\Sph^3\vee\Sph^3\vee\Sph^3\vee\Sph^4\vee\Sph^4$
(see~\cite[Ex.~8.15]{bu-pa02} and~\cite{gr-th04}). Like in the
previous example, the subgroup $\G$ is the diagonal subtorus in
$(\C^*)^3$. By Theorem~\ref{coxth}, $X_\Sigma=U(\Sigma)/\G$, a
quasiprojective variety obtained by removing three points from~$\C
P^2$. This in non-compact, and cannot be identified with $\mathcal
Z(\Sigma)/\K$.
\end{exam}

\section{Normal fans}\label{norfa}
The next step in our study of the Kempf--Ness set for torus
actions on quasiaffine varieties $U(\Sigma)$ would be to obtain an
explicit description like the one given by~\eqref{eqcns} in the
affine case. Although we do not know of such a description in
general, it does exist in the particular case when $\Sigma$ is the
normal fan of a simple polytope.

Let $M_\R=(N_\R)^*$ be the dual vector space. Assume we are given
primitive vectors $a_1,\ldots,a_m\in N$ and integer numbers
$b_1,\ldots,b_m\in\Z$, and consider the set
\begin{equation}\label{ptope}
  P=\{x\in M_\R\colon\langle a_i,x\rangle+b_i\ge0,\quad
  i=1,\ldots,m\}.
\end{equation}
We further assume that $P$ is bounded, the affine hull of $P$ is
the whole $M_\R$, and the intersection of $P$ with every
hyperplane determined by the equation $(a_i,x)+b_i=0$ spans an
affine subspace of dimension $n-1$, for $i=1,\ldots,m$ (or,
equivalently, none of the inequalities can be removed without
enlarging~$P$). This means that $P$ is a \emph{convex polytope}
with exactly $m$ \emph{facets}. (In general, the set $P$ is always
convex, but it may be unbounded, not of the full dimension, or
there may be redundant inequalities.) By introducing a Euclidean
metric in $N_\R$ we may think of $a_i$ as the inward pointing
normal vector to the corresponding facet $F_i$ of $P$, \
$i=1,\ldots,m$. Given a face $Q\subset P$ we say that $a_i$ is
\emph{normal} to $Q$ if $Q\subset F_i$. If $Q$ is an
$q$-dimensional face, then the set of all its normal vectors
$\{a_{i_1},\ldots,a_{i_k}\}$ spans an $(n-q)$-dimensional cone
$\sigma_Q$. The collection of cones $\{\sigma_Q\colon Q\text{ a
face of }P\}$ is a complete fan in $N$, which we denote $\Sigma_P$
and refer to as the \emph{normal fan} of $P$. The normal fan is
simplicial if and only if the polytope $P$ is \emph{simple}, that
is, there exactly $n$ facets meeting at each of its vertices. In
this case the cones of $\Sigma_P$ are generated by subsets
$\{a_{i_1},\ldots,a_{i_k}\}$ such that the intersection
$F_{i_1}\cap\ldots\cap F_{i_k}$ of the corresponding facets is
non-empty.

The Kempf--Ness sets (or the moment angle complexes) $\mathcal
Z(\Sigma_P)$ corresponding to normal fans of simple polytopes
admit a very transparent interpretation as \emph{complete
intersections} of \emph{real algebraic quadrics}, as described
in~\cite{b-p-r??} (these complete intersections of quadrics were
also studied in~\cite{bo-me??}). We give this construction below.

In the rest of this section we assume that $P$ is a simple
polytope, and therefore, $\Sigma_P$ is a simplicial fan. We may
specify $P$ by a matrix inequality $A_Px+b_P\ge0$, where $A_P$ is
the $m\times n$ matrix of row vectors $a_i$, and $b_P$ is the
column vector of scalars $b_i$. The linear transformation
$M_\R\to\R^m$ defined by the matrix $A_P$ is exactly the one
obtained from the map $\mathbb T^m\to\T$ in~\eqref{kgrou} by
applying $\mathop{\mathrm{Hom}}_\Z(\:\cdot\:,\Sph^1)\otimes_\Z\R$.
Since the points of $P$ are specified by the constraint
$A_Px+b_P\ge0$, the formula $i_P(x)=A_Px+b_P$ defines an affine
injection
\begin{equation}\label{defip}
i_P\colon M_\R\longrightarrow\R^m,
\end{equation}
which embeds $P$ into the positive cone $\R^m_\ge=\{y\in\R^m\colon
y_i\ge0\}$.

Now define the space $\mathcal Z_P$ by a pullback diagram
\begin{equation}\label{cdiz}
\begin{CD}
  \mathcal Z_P @>i_Z>>\C^m\\
  @V\varrho_P VV\hspace{-0.2em} @VV\varrho V @.\\
  P @>i_P>> \R^m
\end{CD}
\end{equation}
where $\varrho(z_1,\ldots,z_m)$ is given by
$(|z_1|^2,\ldots,|z_m|^2)$. The vertical maps above are
projections onto the quotients by the $\mathbb T^m$-actions, and
$i_Z$ is a $\mathbb T^m$-\-equi\-va\-ri\-ant embedding.

\begin{prop}{\rm(a)} We have $\mathcal Z_P\subset U(\Sigma_P)$.

{\rm(b)} There is a $\mathbb T^m$-equivariant homeomorphism
$\mathcal Z_P\cong\mathcal Z(\Sigma_P)$.
\end{prop}
\begin{proof}
Assume $z\in\mathcal Z_P\subset\C^m$ and let $\omega(z)$ be the
set of zero coordinates of $z$. Since the facet $F_i$ of $P$ is
the intersection of $P$ with the hyperplane $(a_i,x)+b_i=0$, the
point $\varrho_P(z)$ belongs to the intersection
$\cap_{i\in\omega(z)}F_i$, which is thereby non-empty. Therefore,
the vectors $\{a_i\colon i\in\omega(z)\}$ span a cone of
$\Sigma_P$. Thus, $\omega(z)$ is a $g$-subset and $z\in
U(\Sigma_P)$, proving~(a).

To prove~(b) we look more closely on the construction of the
identification space from the proof of Theorem~\ref{zspro} in the
case when $\Sigma$ is a normal fan. Then the space $C(\Sigma_P)$
may be identified with $P$, and $C(\sigma)$ is the face
$\cap_{i\in g(\sigma)}F_i$ of $P$. The Kempf--Ness set $\mathcal
Z(\Sigma_P)$ is therefore identified with
\begin{equation}\label{idspa}
  (\mathbb T^m\times P)/{\sim}.
\end{equation}
Now we notice that if we replace $P$ by the positive cone
$\R^m_\ge$ (with the obvious face structure) in the above
identification space, we obtain $(\mathbb
T^m\times\R^m_\ge)/{\sim}=\C^m$. Since the map $i_P$
from~\eqref{cdiz} respects facial codimension, the pullback space
$\mathcal Z_P$ can be also identified with~\eqref{idspa}, thus
proving~(b).
\end{proof}

Choosing a basis for $\mathop{\mathrm{coker}}A_P$ we obtain a
$(m-n)\times m$-matrix $C$ so that the resulting short exact
sequence
\begin{equation}\label{pexse}
  0\longrightarrow M_\R\stackrel{A_P}{\longrightarrow}\R^m
  \stackrel{C}{\longrightarrow}\R^{m-n}\longrightarrow 0,
\end{equation}
is the one obtained from~\eqref{kgrou} by applying
$\mathop{\mathrm{Hom}}_\Z(\:\cdot\:,\Sph^1)\otimes_\Z\R$.

We may assume that the first $n$ normal vectors $a_1,\ldots,a_n$
span a cone of $\Sigma_P$ (equivalently, the corresponding facets
of $P$ meet at a vertex), and take these vectors as a basis of
$M_\R$. In this basis, the first $n$ rows of the matrix $(a_{ij})$
of $A_P$ form a unit $n\times n$-matrix, and we may take
\begin{equation}\label{cmatr}
  C=\begin{pmatrix}
    -a_{n+1,1} & \ldots & -a_{n+1,n} & 1     & 0      & \ldots & 0\\
    -a_{n+2,1} & \ldots & -a_{n+2,n} &      0 & 1     & \ldots & 0\\
    \vdots    & \ddots & \vdots    & \vdots &\vdots&\ddots&\vdots\\
    -a_{m,1}   & \ldots & -a_{m,n}   & 0 & 0  & \ldots & 1
  \end{pmatrix}.
\end{equation}
Then diagram~\eqref{cdiz} implies that $i_Z$ embeds $\mathcal Z_P$
in $\C^m$ as the set of solutions of the $m-n$ real quadratic
equations
\begin{equation}\label{zpqua}
  \sum_{k=1}^mc_{jk}\left(|z_k|^2-b_k\right)=0,\;\text{ for }
  1\leq j\leq m-n
\end{equation}
where $C=(c_{jk})$ is given by~\eqref{cmatr}. This intersection of
real quadrics is non-degenerate~\cite[Lemma~3.2]{b-p-r??} (the
normal vectors are linearly independent at each point), and
therefore, $\mathcal Z_P\subset\R^{2m}$ is a smooth submanifold
with trivial normal bundle.

\section{Projective toric varieties and moment maps}\label{prvar}
In the notation of Section~\ref{knsav}, let
$f_v=(dF_v)_e\colon\mathfrak g\to\R$. It maps $\gamma\in\mathfrak
g$ to $\mathop{\mathrm{Re}}\langle\gamma v,v\rangle$,
see~\eqref{eqcns}. We may consider $f_v$ as an element of the dual
Lie algebra $\mathfrak g^*$. As $\G$ is reductive, we have
$\mathfrak g=\mathfrak k\oplus i\mathfrak k$. Since $\K$ is norm
preserving, $f_v$ vanishes on $\mathfrak k$, so we consider $f_v$
as an element of $i\mathfrak k^*\cong\mathfrak k^*$. Varying $v\in
V$ we get the \emph{moment map} $\mu\colon V\to\mathfrak k^*$,
which sends $v\in V$, $\kappa\in\mathfrak k$ to $\langle i\kappa
v,v\rangle$. The Kempf--Ness set is the set of zeroes of $\mu$:
\begin{equation}\label{muinv}
  \kn=\mu^{-1}(0).
\end{equation}

This description does not apply to the case of algebraic torus
actions on $U(\Sigma)$ considered in the two previous sections: as
is seen from simple examples below, the set $\mu^{-1}(0)=\{z\in
\C^m\colon\langle\kappa z,z\rangle=0\text{ for all
}\kappa\in\mathfrak k\}$ consists only of the origin in this case.
Nevertheless, in this section we show that a description of the
toric Kempf--Ness set $\mathcal Z(\Sigma)$ similar
to~\eqref{muinv} exists in the case when $\Sigma$ is a normal fan,
thereby extending the analogy with Kempf--Ness sets for affine
varieties even further.

As explained in~\cite{dani78} or~\cite[\S5.1]{bu-pa02}, the toric
variety $X_\Sigma$ is projective exactly when $\Sigma$ arises as
the normal fan of a convex polytope. In fact, the set of integers
$\{b_1,\ldots,b_m\}$ from~\eqref{ptope} determines an \emph{ample}
divisor on $X_{\Sigma_P}$, therefore providing a projective
embedding. Note that the vertices of $P$ are not necessarily
lattice points in $M$ (as they may have rational coordinates), but
this can be remedied by simultaneously multiplying
$b_1,\ldots,b_m$ by an integer number; this corresponds to the
passage from an ample divisor to a \emph{very ample} one.

Assume now that $\Sigma_P$ is a regular fan, therefore,
$X_{\Sigma_P}$ is a smooth projective variety. This implies that
$X_{\Sigma_P}$ is K\"ahler, and therefore, a symplectic manifold.
There is the following symplectic version of the construction from
Section~\ref{altoa}.

Let $(W,\omega)$ be a symplectic manifold with a $\K$-action that
preserves the symplectic form $\omega$. For every
$\kappa\in\mathfrak k$ we denote by $\xi_\kappa$ the corresponding
$\K$-invariant vector field on $W$. The $\K$-action is said to be
\emph{Hamiltonian} if the 1-form $\omega(\:\cdot\:,\xi_\kappa)$ is
exact for every $\kappa\in\mathfrak k$, that is, there is a
function $H_\kappa$ on $W$ such that
\[
  \omega(\xi,\xi_\kappa)=dH_\kappa(\xi)=\xi(H_\kappa)
\]
for every vector field $\xi$ on $W$. Under this assumption, the
\emph{moment map}
\[
  \mu\colon W\to\mathfrak k^*,\qquad (x,\kappa)\mapsto H_\kappa(x)
\]
is defined.

\begin{exam}
1. Basic example is given by $W=\C^m$ with the symplectic form
$\omega=2\sum_{k=1}^mdx_k\wedge dy_k$, where $z_k=x_k+iy_k$. The
coordinatewise action of $\mathbb T^m$ is Hamiltonian with the
moment map $\mu\colon\C^m\to\R^m$ given by
$\mu(z_1,\ldots,z_m)=(|z_1|^2,\ldots,|z_m|^2)$ (we identify the
dual Lie algebra of $\mathbb T^m$ with $\R^m$).

2. Now let $\Sigma$ be a regular fan, and $\K$ be the subgroup of
$\mathbb T^m$ defined by~\eqref{kgrou}. We can restrict the
previous example to the $\K$-action on the invariant subvariety
$U(\Sigma)\subset\C^m$. The corresponding moment map is then
defined by the composition
\begin{equation}\label{tmoma}
  \mu_\Sigma\colon\C^m\longrightarrow\R^m\longrightarrow \mathfrak k^*.
\end{equation}
A choice of an isomorphism $\mathfrak k\cong\R^{m-n}$ allows to
identify the map $\R^m\to\mathfrak k^*$ with the linear
transformation given by matrix~\eqref{cmatr}, see~\eqref{pexse}.
\end{exam}

A direct comparison with \eqref{muinv} prompts us to relate the
level set $\mu_\Sigma^{-1}(0)$ of moment map~\eqref{tmoma} with
the toric Kempf--Ness set $\mathcal Z(\Sigma_P)$ for the
$\G$-action on $U(\Sigma_P)$. However, this analogy is not that
straightforward: the set $\mu^{-1}_{\Sigma}(0)=\{z\in
\C^m\colon\langle\kappa z,z\rangle=0\text{ for all
}\kappa\in\mathfrak k\}$ is given by the equations
$\sum_{k=1}^mc_{jk}|z_k|^2=0$, \ $1\leq j\leq m-n$, which have
only zero solution. (Indeed, as the intersection of $\R^m_\ge$
with the affine $n$-plane $i_P(M_\R)=A_P(M_\R)+b_P$ is bounded,
its intersection with the plane $A_P(M_\R)$ consists only of the
origin.) On the other hand, by comparing~\eqref{tmoma}
with~\eqref{zpqua}, we obtain

\begin{prop}
Let $\Sigma_P$ be the normal fan of a simple polytope given
by~\eqref{ptope}, and~\eqref{tmoma} the corresponding moment map.
Then the toric Kempf--Ness set $\mathcal Z(\Sigma_P)$ for the
$\G$-action on $U(\Sigma_P)$ is given by
\[
  \mathcal Z(\Sigma_P)\cong\mu^{-1}_{\Sigma_P}(Cb_P).
\]
\end{prop}

In other words, the difference with the affine situation is that
we have to take $Cb_P$ instead of 0 as the value of the moment
map. The reason is that $Cb_P$ is a \emph{regular value} of $\mu$,
unlike 0.

By making a perturbation $b_i\mapsto b_i+\varepsilon_i$ of the
values $b_i$ in~\eqref{ptope} while keeping the vectors $a_i$
unchanged for $1\le i\le m$, we obtain another convex set
$P(\varepsilon)$ determined by~\eqref{ptope}. Provided that the
perturbation is small, the set $P(\varepsilon)$ is still a simple
convex polytope of the same \emph{combinatorial type} as~$P$. Then
the normal fans of $P$ and $P(\varepsilon)$ are the same, and the
manifolds $\mathcal Z_P$ and $\mathcal Z_{P(\varepsilon)}$ defined
by~\eqref{zpqua} are $\mathbb T^m$-equivariantly homeomorphic.
Moreover, $Cb_{P(\varepsilon)}$, considered as an element of
$\mathfrak
k^*=\mathop{\mathrm{Hom}}_\Z(\K,\Sph^1)\otimes_\Z\R\cong
H^2(X_{\Sigma_P};\R)$, belongs to the \emph{K\"ahler cone} of the
toric variety $X_{\Sigma_P}$~\cite[\S4]{cox97}. In the case of
normal fans the following version of our Theorem~\ref{zspro}~(a)
is known in toric geometry:

\begin{theo}[{cf.~\cite[Th.~4.1]{cox97}}]
Let $X_\Sigma$ be a projective simplicial toric variety and assume
that $c\in H^2(X_{\Sigma};\R)$ is in the K\"ahler cone. Then
$\mu_\Sigma^{-1}(c)\subset U(\Sigma)$, and the natural map
\[
  \mu_\Sigma^{-1}(c)/\K\to U(\Sigma)/\G=X_\Sigma
\]
is a diffeomorphism.
\end{theo}

This statement is the essence of the construction of smooth
projective toric varieties via \emph{symplectic reduction}. The
submanifold $\mu_\Sigma^{-1}(c)\subset\C^m$ may fail to be
symplectic because the restriction of the standard symplectic form
$\omega$ on $\C^m$ to $\mu_\Sigma^{-1}(c)$ may be degenerate.
However, the restriction of $\omega$ descends to the quotient
$\mu_\Sigma^{-1}(c)/\K$ as a symplectic form.

\begin{exam}
Let $P=\Delta^n$ be the \emph{standard simplex} defined by $n+1$
inequalities $\langle e_i,x\rangle\ge0$, \ $i=1,\ldots,n$, and
$\langle-e_1-\ldots-e_n,x\rangle+1\ge0$ in $M_\R$ (here
$e_1,\ldots,e_n$ is a chosen basis which we use to identify $N_\R$
with $\R^n$). The cones of the corresponding normal fan $\Sigma$
are generated by the proper subsets of the set of vectors
$\{e_1,\ldots,e_n,-e_1-\ldots-e_n\}$. The groups $\G\cong\C^*$ and
$\K\cong\Sph^1$ are the diagonal subgroups in $(\C^*)^{n+1}$ and
$\mathbb T^{n+1}$ respectively, while
$U(\Sigma)=\C^{n+1}\setminus\{0\}$. The $(n+1)\times n$-matrix
$A_P=(a_{ij})$ has $a_{ij}=\delta_{ij}$ for $1\le i,j\le n$ and
$a_{n+1,j}=-1$ for $1\le j\le n$. The matrix $C$~\eqref{cmatr} is
just one row of units. Moment map~\eqref{tmoma} is given by
$\mu_\Sigma(z_1,\ldots,z_{n+1})=|z_1|^2+\ldots+|z_{n+1}|^2$. Since
$Cb_P=1$, the Kempf--Ness set $\mathcal Z_P=\mu^{-1}_\Sigma(1)$ is
the unit sphere $\Sph^{2n+1}\subset\C^{n+1}$, and
$X_\Sigma=(\C^{n+1}\setminus\{0\})/\G=\Sph^{2n+1}/\K$ is the
complex projective space $\C P^n$.
\end{exam}

In the next section we consider a more complicated example, while
here we finish with an open question.

\begin{prob}
As is known (see e.g.~\cite[Ch.~5]{bu-pa02}), there are many
complete regular fans $\Sigma$ which cannot be realised as normal
fans of convex polytopes. The corresponding toric varieties
$X_\Sigma$ are not projective (although being non-singular). In
this case the toric Kempf--Ness set $\mathcal Z(\Sigma)$ is still
defined (see Section~\ref{altoa}). However, the rest of the
analysis of the last two sections does not apply here; in
particular, we don't have a description of $\mathcal Z(\Sigma)$ as
in~\eqref{zpqua}. Can one still describe $\mathcal Z(\Sigma)$ as a
complete intersection of real quadratic (or higher order)
hypersurfaces?
\end{prob}

\section{Cohomology of toric Kempf--Ness sets}\label{cohom}
Here we use the results of~\cite{bu-pa02} and~\cite{pano05} on
moment-angle complexes to describe the integer cohomology rings of
toric Kempf--Ness sets. As we shall see from an example below, the
topology of $\mathcal Z(\Sigma)$ may be quite complicated even for
simple fans.

Given an abstract simplicial complex $\mathcal K$ on the set
$[m]=\{1,\ldots,m\}$, the \emph{face ring} (or the
\emph{Stanley--Reisner ring}) $\Z[\mathcal K]$ is defined as the
following quotient of the polynomial ring on $m$ generators:
\[
  \Z[\mathcal K]=\Z[v_1,\ldots,v_m]/(v_{i_1}\cdots v_{i_k}\colon\{i_1,\ldots,i_k\}\text{ is not a simplex of
  }\mathcal K).
\]
We introduce a grading by setting $\deg v_i=2$, \ $i=1,\ldots,m$.
As $\Z[\mathcal K]$ may be thought of as a
$\Z[v_1,\ldots,v_m]$-module via the projection map, the bigraded
\emph{$\Tor$-modules}
$\Tor_{\Z[v_1,\ldots,v_m]}^{-i,2j}(\Z[\mathcal K],\Z)$ are
defined, see~\cite{stan96}. They can be calculated, for example,
using the \emph{Koszul resolution} of the trivial
$\Z[v_1,\ldots,v_m]$-module $\Z$. This also endows
$\Tor^*_{\Z[v_1,\ldots,v_m]}(\Z[\mathcal K],\Z)$ with a graded
commutative algebra structure (the grading is by the total
degree), see details in~\cite[Ch.~7]{bu-pa02}.

\begin{theo}[{cf.~\cite[Th.~7.6,~7.7]{bu-pa02},~\cite[Th.~4.7]{pano05}}]\label{zkcoh}
For every simplicial fan $\Sigma$ there are algebra isomorphisms
\[
  H^*\bigl(\mathcal Z(\Sigma);\Z\bigr)
  \cong\Tor^*_{\Z[v_1,\ldots,v_m]}\bigl(\Z[\mathcal K_\Sigma],\Z\bigr)\cong
  H\bigl[\Lambda[u_1,\ldots,u_m]\otimes\Z[\mathcal K_\Sigma],d\bigr],
\]
where the latter denotes the cohomology of a differential graded
algebra with $\deg u_i=1$, $\deg v_i=2$, $du_i=v_i$, $dv_i=0$ for
$1\le i\le m$.
\end{theo}

Given a subset $I\subseteq[m]$, denote by $\mathcal K(I)$ the
corresponding \emph{full subcomplex} of $\mathcal K$, or the
restriction of $\mathcal K$ to $I$. We also denote by $\widetilde
H^i(\mathcal K(I))$ the $i$th reduced simplicial cohomology group
of $\mathcal K(I)$ with integer coefficients. A theorem due to
Hochster~\cite{hoch77} expresses the $\Tor$-modules
$\Tor_{\Z[v_1,\ldots,v_m]}^{-i,2j}(\Z[\mathcal K],\Z)$ in terms of
full subcomplexes of $\mathcal K$, which leads to the following
description of the cohomology of $\mathcal Z(\Sigma)$.

\begin{theo}[{cf.~\cite[Cor.~5.2]{pano05}}]\label{hochs}
We have
\[
  H^k\bigl(\mathcal Z(\Sigma)\bigr)\cong\bigoplus_{I\subseteq[m]}
  \widetilde H^{k-|I|-1}\bigl(\mathcal K_\Sigma(I)\bigr).
\]
\end{theo}
\noindent There is also a description of the product in
$H^*(\mathcal Z(\Sigma))$ in terms of full subcomplexes of
$\mathcal K_\Sigma$, cf.~\cite[Th.~5.1]{pano05}.

\begin{exam}
Let $P$ be the simple polytope obtained by cutting two
non-adjacent edges off a cube in $M_\R\cong\R^3$, as shown on
Fig.~\ref{figur}. We may specify such a polytope by 8 inequalities
\begin{gather*}
  x\ge0,\quad y\ge0,\quad
  z\ge0,\quad-x+3\ge0,\quad-y+3\ge0,\quad-z+3\ge0,\\
  -x+y+2\ge0,\quad -y-z+5\ge0,
\end{gather*}
and it has 8 facets $F_1,\ldots,F_8$, numbered as on the picture.
\begin{figure}[h]
\begin{center}
\begin{picture}(75,80)
  \put(10,15){\line(2,-1){10}}
  \put(10,15){\line(0,1){40}}
  \put(20,50){\line(-2,1){10}}
  \put(20,50){\line(0,-1){40}}
  \put(20,50){\line(1,0){20}}
  \put(40,50){\line(1,-1){5}}
  \put(45,45){\line(0,-1){35}}
  \put(45,10){\line(-1,0){25}}
  \put(45,10){\line(1,1){20}}
  \put(65,30){\line(0,1){35}}
  \put(65,65){\line(-1,-1){20}}
  \put(40,50){\line(1,1){20}}
  \put(60,70){\line(1,-1){5}}
  \put(60,70){\line(-1,0){35}}
  \put(25,70){\line(-1,-1){15}}
  \multiput(25,30)(5.1,0){8}{\line(1,0){3}}
  \multiput(25,30)(0,5.1){8}{\line(0,1){3}}
  \multiput(25,30)(-5.2,-5.2){3}{\line(-1,-1){3.6}}
  \put(25,26.5){0}
  \put(10,15){\vector(-1,-1){10}}
  \put(25,70){\vector(0,1){10}}
  \put(65,30){\vector(1,0){10}}
  \put(-2,6){$x$}
  \put(75,27){$y$}
  \put(26,78){$z$}
  \put(14,30){$F_7$}
  \put(49,55){$F_8$}
  \put(30,58){$F_6$}
  \put(35,6.5){\vector(0,1){3}}
  \put(34,3){$F_3$}
  \put(54,35){$F_5$}
  \put(13,66){\vector(1,-1){3.6}}
  \put(10,67){$F_2$}
  \put(30,32){$F_4$}
  \put(50,75){\vector(0,-1){4}}
  \put(48,76){$F_1$}
\end{picture}
\caption{ } \label{figur}
\end{center}
\end{figure}
The 1-dimensional cones of the corresponding normal fan $\Sigma_P$
are spanned by the following primitive vectors:
\begin{gather*}
  a_1=e_1,\quad a_2=e_2,\quad a_3=e_3,\quad
  a_4=-e_1,\quad a_5=-e_2,\quad a_6=-e_3,\\
  a_7=-e_1+e_2,\quad a_8=-e_2-e_3.
\end{gather*}
Toric variety $X_{\Sigma_P}$ is obtained by blowing up the product
$\C P^1\times\C P^1\times\C P^1$ (corresponding to the cube) at
two complex 1-dimensional subvarieties
$\{\infty\}\times\{0\}\times\C P^1$ and $\C
P^1\times\{\infty\}\times\{\infty\}$. Matrix~\eqref{cmatr} is
given by
\[
  C=\begin{pmatrix}
    1&0&0&1&0&0&0&0\\
    0&1&0&0&1&0&0&0\\
    0&0&1&0&0&1&0&0\\
    1&-1&0&0&0&0&1&0\\
    0&1&1&0&0&0&0&1
  \end{pmatrix}.
\]
Its transpose determines the inclusion $\G\hookrightarrow(\C^*)^8$
(or $\K\hookrightarrow T^8$), and we have
$X_{\Sigma_P}=U(\Sigma_P)/\G=\mathcal Z(\Sigma_P)/\K$ by
Theorem~\ref{zspro}. The toric Kempf--Ness set $\mathcal
Z(\Sigma_P)\cong\mathcal Z_P$~\eqref{zpqua} is defined by 5 real
quadratic equations:
\begin{gather*}
|z_1|^2+|z_4|^2-3=0,\quad |z_2|^2+|z_5|^2-3=0,\quad
|z_3|^2+|z_6|^2-3=0,\\
|z_1|^2-|z_2|^2+|z_7|^2-2=0,\quad |z_2|^2+|z_3|^2+|z_8|^2-5=0.
\end{gather*}

The dual triangulation $\mathcal K_{\Sigma}$ is obtained from the
boundary of an octahedron by applying two stellar subdivisions at
non-adjacent edges~\cite{ma-pa06}. The face ring is
\[
  \Z[\mathcal K_\Sigma]\!=\!\Z[v_1,\ldots,v_8]/(v_1v_4,v_1v_7,v_2v_4,v_2v_5,
  v_2v_8,v_3v_6,v_3v_8,v_5v_6,v_5v_7,v_7v_8).
\]
According to Theorem~\ref{hochs}, the group $H^3(\mathcal Z_P)$
has a generator for every pair of vertices of $\mathcal K_\Sigma$
not joined by an edge (equivalently, for every pair of
non-adjacent facets of $P$). Therefore, $H^3(\mathcal
Z_P)\cong\Z^{10}$, and the generators are represented by the
following 3-cocycles in the differential graded algebra from
Theorem~\ref{zkcoh}:
\[
  u_1v_4,\;u_1v_7,\;u_2v_4,\;u_2v_5,\;
  u_2v_8,\;u_3v_6,\;u_3v_8,\;u_5v_6,\;u_5v_7,\;u_7v_8.
\]

Using Theorem~\ref{hochs} again, we see that only reduced
0-cohomology of 3-vertex full subcomplexes of $\mathcal K_\Sigma$
may contribute to $H^4(\mathcal Z_P)$. There are two types of
disconnected simplicial complexes on 3 vertices: ``3 disjoint
points'' and ``an edge and a point''. $\mathcal K_\Sigma$ contains
no full subcomplexes of the first type and 16 subcomplexes of the
second type. The corresponding 4-cocycles in the differential
graded algebra $\Lambda[u_1,\ldots,u_m]\otimes\Z[\mathcal
K_\Sigma]$ are
\begin{gather*}
  u_4u_7v_1,\;u_4u_5v_2,\;u_4u_8v_2,\;u_5u_8v_2,\;
  u_6u_8v_3,\;u_1u_2v_4,\;u_2u_6v_5,\;u_2u_7v_5,\\
  u_6u_7v_5,\;u_3u_5v_6,\;u_1u_5v_7,\;u_1u_8v_7,\;
  u_5u_8v_7,\;u_2u_3v_8,\;u_2u_7v_8,\;u_3u_7v_8.
\end{gather*}
Therefore, $H^4(\mathcal Z_P)\cong\Z^{16}$.

The 5th cohomology group of $\mathcal Z_P$ is the sum of the 1st
cohomology of 3-vertex full subcomplexes of $\mathcal K_\Sigma$
and the reduced 0-cohomology of 4-vertex full subcomplexes. A
3-vertex full subcomplex of $\mathcal K_\Sigma$ may have non-zero
1st cohomology group only if the corresponding 3 facets of $P$
form a ``belt'', that is, are pairwise adjacent but do not share a
common vertex. As there are no such 3-facet belts in $P$, only
reduced 0-cohomology of 4-vertex subcomplexes contributes in
$H^5(\mathcal Z_P)$. The corresponding 5-cocycles are
\begin{gather*}
  u_1u_5u_8v_7,\;u_2u_3u_7v_8,\;u_4u_5u_8v_2,\;u_2u_6u_7v_5,\;
  u_2u_7u_5v_8-u_2u_7u_8v_5
\end{gather*}
(note that the last cocycle cannot be represented by a monomial).
Therefore, $H^5(\mathcal Z_P)\cong\Z^5$. Due to Poincar\'e
duality, this completely determines the Betti vector
$(1,0,0,10,16,5,5,16,10,0,0,1)$ of the 11-dimensional manifold
$\mathcal Z_P$. The generators of the sixth cohomology group,
$H^6(\mathcal Z_P)\cong\Z^5$, correspond to the 4-facet belts in
$P$, and the corresponding 6-cocycles are
\begin{gather*}
  u_2u_3v_4v_6,\;u_1u_5v_4v_6,\;u_1u_3v_6v_7,\;u_1u_3v_4v_8,\;
  u_1u_3v_4v_6.
\end{gather*}
These are the Poincar\'e duals to the 5-cocycles. The fundamental
class of $\mathcal Z_P$ is represented (up to a sign) by the
cocycle $u_4u_5u_6u_7u_8v_1v_2v_3$, or by any cocycle of the form
$u_{\sigma(4)}u_{\sigma(5)}u_{\sigma(6)}u_{\sigma(7)}
u_{\sigma(8)}v_{\sigma(1)}v_{\sigma(2)}v_{\sigma(3)}$ where
$\sigma\in S_8$ is a permutation such that the facets
$F_{\sigma(1)}$, $F_{\sigma(2)}$ and $F_{\sigma(3)}$ share a
common vertex.

The multiplicative structure in $H^*(\mathcal Z_P)$ can be easily
retrieved from this description. For example, we have identities
\begin{gather*}
  [u_1v_4]\cdot[u_1v_7]=0,\quad[u_1v_7]\cdot[u_2v_4]=0,\quad
  [u_1v_4]\cdot[u_3v_6]=[u_1u_3v_4v_6],\\
  [u_2v_4]\cdot[u_3v_6]\cdot[u_1u_5u_8v_7]=
  [u_1u_2u_3u_5u_8v_4v_6v_7],\quad
  \text{etc.}
\end{gather*}

Yet another interesting feature of the manifold $\mathcal Z_P$ of
this Example is the existence of non-trivial Massey products in
$H^*(\mathcal Z_P)$~\cite{bask03}. Consider 3 cocycles $a=u_1v_4$,
\ $b=u_2v_5$, \ $c=u_3v_6$ representing cohomology classes
$\alpha,\beta,\gamma\in H^3(\mathcal Z_P)$. Since $\alpha\beta=0$
and $\beta\gamma=0$, a triple Massey product
$\langle\alpha,\beta,\gamma\rangle$ is defined. It consists of the
cohomology classes in $H^8(\mathcal Z_P)$ represented by the
cocycles of the form $af+ec$ for all choices of $e$ and $f$ such
that $ab=de$ and $bc=df$ (here $d$ denotes the differential; as
there are may be many choices of $e$ and $f$, the Massey product
is a multivalued operation in general). The Massey product is said
to be \emph{trivial} if it contains zero. In our case we may take
$e=u_1u_2u_5v_4$ and $f=0$, so $\langle\alpha,\beta,\gamma\rangle$
contains a non-zero cohomology class $[u_1u_2u_5u_3v_4v_6]\in
H^8(\mathcal Z_P)$. Moreover, $\langle\alpha,\beta,\gamma\rangle$
is non-trivial, cf.~\cite[Ex.~5.7]{pano05}. This implies that
$\mathcal Z_P$ is a \emph{non-formal} manifold. A detailed study
of Massey products in the cohomology of moment-angle complexes is
undertaken in~\cite{de-su05}.
\end{exam}

\end{document}